\documentclass[11pt]{article}
\usepackage{mathrsfs}
\usepackage{amssymb}
\usepackage{amsfonts}
\usepackage{color,xcolor}
\usepackage{graphicx}
\usepackage{geometry}
\usepackage{manfnt}
\RequirePackage[colorlinks,citecolor=blue,urlcolor=blue]{hyperref}
\newtheorem{theorem}{Theorem}[section]

\newtheorem{lemma}[theorem]{Lemma}

\newtheorem{remark}{Remark}[section]

\def\[{{\Big[}}\def\]{{\Big]}}\def\<{{\langle}}\def\>{{\rangle}}\def\({{\Big(}}
\def\){{\Big)}}

\def\min{{\mathord{{\rm min}}}}

\def\={&\!\!=\!\!&}

\def\geq{\geqslant}\def\leq{\leqslant}

\begin{document}
\title{\bf
Eulerian collinear configuration for 3-body problem}
\author{Liang Ding$^{1}$, Juan Manuel S\'{a}nchez-Cerritos$^{2}$ and Jinlong Wei$^{3\,\ast}$
\\ {\small \it $^1$School of Data Science and Information Engineering, Guizhou} \\  {\small \it Minzu University, Guiyang, 550025, P.R. China}\\
{\small \tt ding2016liang@126.com}\\  {\small \it $^2$College of Mathematics and Statistics, Chongqing Technology} \\  {\small \it and Business University, Chongqing 400067, P.R. China}\\
 {\small \tt sanchezj011@outlook.com}\\  {\small \it $^3$School of Statistics and Mathematics, Zhongnan University} \\  {\small \it of Economics and Law, Wuhan, 430073, P.R. China}\\
 {\small \tt weijinlong.hust@gmail.com}}

\date{}
 \maketitle
\noindent{\hrulefill} \vskip1mm\noindent
 {\bf Abstract}
For 3-body problem with any given masses $m_1, \,m_2,\,m_3>0$, there exist only Eulerian collinear central configuration and Lagrangian equilateral-triangle central configuration, and in this paper, for planar 3-body problem,
we prove that there exists another non-collision trajectory $q$, which is not the variational
minimizer of the Lagrangian action on the loop space $\overline{\Lambda_1}$,
is also an Eulerian collinear central configuration at any instant. Moreover, we do not need the restriction condition on the winding number $deg(q_i-q_j)\neq0 \,(i\neq j)$.

 \vskip1mm\noindent{\hrulefill}

\section{Introduction} \label{sec1}
\setcounter{equation}{0}
For the motion of 3-body problem, it is well-known that Newton's second law and Newton's universal gravitation law yield the following equation
\begin{eqnarray}\label{1.1}
m_i\ddot{q}_i+\sum_{{j\neq i\atop{{1\leq j\leq 3}}}}m_im_j\frac{q_i-q_j}{|q_j-q_i|^3}=0, \quad q_{i}\in\mathbb{R}^{k}, \quad k=2 \,\,or\,\,3, \,\,i=1,2,3,
\end{eqnarray}
and when $k=2$, we call it planar 3-body problem; when $k=3$, we call it spatial 3-body problem.
In this paper, for planar 3-body problem (\ref{1.1}), we attempt to study the non-collision solution with central configuration characteristic, and firstly, we introduce a definition on central configuration.

\smallskip \noindent
\textbf{Definition 1.1.}(\cite{Wintner1947}) At a given instant $t=t_{0}$, for $N$ mass points $m_i$ with positions $q_{i}(t_{0}) \in \mathbb{R}^{k}$
$(k=2\,\,or\,\,3,\,i=1,2,\ldots,N)$, a configuration $q(t_{0})=(q_1(t_{0}),q_2(t_{0}),\ldots,q_{N}(t_{0}))^{T}\in X\setminus \Delta$ is called a central configuration if there exists a constant $\lambda\in \mathbb{R}$ such that
\begin{eqnarray*}
\left\{\begin{array}{ll}
\sum\limits_{{j\neq i\atop{{1\leq j\leq N}}}}
\frac{m_j m_i}{|q_j(t_{0})-q_i(t_{0})|^{3}}(q_j(t_{0})-q_i(t_{0}))=\lambda m_{i}(q_{i}(t_{0})-c_{0}), \quad  i=1,2,\ldots,N, \\
\ \lambda=\frac{V(q(t_{0}))}{I(q(t_{0}))},
\end{array}\right.
\end{eqnarray*}
where
\begin{eqnarray*}
X\setminus \triangle=\{q(t_{0})=(q_{1}(t_{0}),\ldots,q_{N}(t_{0}))^{T}\in (\mathbb{R}^{2})^{N} \,or \, (\mathbb{R}^{3})^{N}: q_{i}(t_{0})\neq q_{j}(t_{0}), \, when \, i\neq j\}.
\end{eqnarray*}
$c_0$, $V$ and $I$ represent the center of masses, the Newtonian potential and the moment of inertia, respectively, which are given by
\begin{eqnarray*}
c_0=\frac{\sum\limits_{1\leq i\leq N}m_{i}q_{i}(t_{0})}{\sum\limits_{1\leq i\leq N}m_{i}},\, V(q(t_{0}))=-\sum_{1\leq i<j\leq N}\frac{m_i m_j}{|q_i(t_{0})-q_j(t_{0})|}
\end{eqnarray*}
and
\begin{eqnarray*}
I(q(t_{0}))=\sum_{1\leq i\leq N}m_i|q_i(t_{0})-c_{0}|^{2},
\end{eqnarray*}
respectively.

The set of central configurations are invariant
under three classes of transformations: translations, scalings and orthogonal
transformations \cite{Wintner1947}, and the study of central configurations is a very important subject in celestial mechanics with a long and varied history \cite{Moeckel1990,Smale1998}, and a well-known fact is that finding the relative equilibrium solutions of the classical $N$-body problem and the planar central configurations are equivalent. But the problem on the numbers for central configurations is a very difficult topic \cite{Albouy2012,Hampton2006}, and Smale took it as one of the most important 18 mathematical problems (the sixth one) for the 21st century \cite{Smale1998}, one of the reasons is that finding the concrete central configurations is a very difficult work \cite{Fernandes2017,Hampton2005,MacMillan1932,Perez2007}.
For $N=3$, in 1767, Euler \cite{E1767} found the well-known Euler's collinear central configuration (Euler's central configuration for short), and in 1772, Lagrange \cite{Lagrange1772} proved the existence of the famous equilateral-triangle central configuration (also called Lagrange's central configuration for short). In fact, for 3-body problem, there are only two kinds of concrete central configurations: Euler's central configuration and Lagrange's central configuration \cite{Abraham1978}, and in this paper, we study one kind of
important concrete central configurations: Euler's central configuration, i.e. we discuss
(\ref{1.1}) in which $q$ can form an Euler's central configuration at any instant. Before recalling
some elegant works which are relevant for the present paper, we introduce a notation on winding number.

\smallskip\noindent
\textbf{Defiition 1.2.} (\cite{Chern1989}) Let $C=\{x(t), \,t\in[a,b]\}$ be a given oriented closed curve, and $p\in\mathbb{R}^{2}$ be a point not on the curve, then the mapping $\phi: \, C\rightarrow S^{1}$, given by
\begin{eqnarray*}
\phi(x(t))=\frac{x(t)-p}{|x(t)-p|}, \,\,t\in[a,\,b],
\end{eqnarray*}
is defined to be a position mapping of the curve $C$ relative to $p$. When the point on $C$ goes around the
curve once, its image $\phi(x(t))$ will go around a number of times, this number is called the winding number
of the curve $C$ relative to $p$, and we denote it by $deg(x)$.

For 2-body problem, in 1977, Gordon \cite{Gordon1977} proved that the elliptic Keplerian orbit minimize the Lagrangian
action of 2-body problem. For planar 3-body problem, there are elegant works: in 2000, by assuming anti-$T/2$ symmetry condition, Long and Zhang \cite{Long2000} (or \cite{Chenciner1998}) proved that for any given positive masses $m_1$, $m_2$ and $m_3$, the variational
minimizer of the Lagrangian action in $\mathbb{R}^{2}$, is  precisely Lagrange's central configuration at any instant; in 2001, by assuming the winding number $deg \,(q_i-q_j)\neq0\,(i\neq j)$, Zhang and Zhou \cite{Zhang2001} (or \cite{Venturelli2001}) also proved that for any given choice of the three positive masses, the variational
minimizer of the Lagrangian action in $\mathbb{R}^{2}$ is the Lagrange's central configuration at any instant. For more details in this direction, we refer to \cite{Llibre2015,Wintner1947,Zhang2004}.

Note that for 3-body problem, there are only two kinds of concrete central configurations: Euler's central configuration and Lagrange's central configuration. We also note that for planar 3-body problem, in 2004, employing a direct variational method with that assumptions that the winding number $deg(q_i-q_j)\neq0\,(i\neq j)$, and the three bodies are collinear, Zhang and Zhou \cite{Zhang2004} proved that the minimizer of Lagrangian action $f(q)$ on $\overline{\Lambda_1}$, just forms an Euler's central configuration at any instant, where the loop space $\Lambda_1$ is defined as follows:
\begin{eqnarray}\label{1.2}
\Lambda_1&=&\{q=(q_1,q_2,q_3) \,\mid \, \,q_i\in W^{1,2}(\mathbb{R}/T\mathbb{Z},\,\mathbb{R}^{2}), \,q_{3}(t)-q_{1}(t)=\lambda_0 (q_{2}(t)-q_{1}(t)), \nonumber \\ && \sum_{i=1}^{3}m_{i}q_{i}=0, \,deg(q_i-q_j)\neq0 \,\,for \,\,i\neq j\},
\end{eqnarray}
where $m_i>0 \ (i=1,2,3)$ and $T>0$ is any given masses and period respectively, $\lambda_{0}$ satisfies
\begin{eqnarray}\label{1.3}
\frac{m_{3}\lambda_{0}^{-2}+m_2}{m_{3}\lambda_{0}+m_2}-\frac{m_{3}(1-\lambda_{0})^{-2}+m_1}{m_{3}(1-\lambda_{0})+m_1}=0.
\end{eqnarray}
Then an interesting question is proposed:

\textbf{Question.} For planar 3-body problem (\ref{1.1}), can we use the anti-$T/2$ symmetry condition to substitute the winding number condition $deg(q_i-q_j)\neq0\,(i\neq j)$, to obtain new existence of non-collision trajectory $q$ such that the three bodies form an Euler's central configuration at any instant, and $q$ is not the variational
minimizer of the Lagrangian action on $\overline{\Lambda_1}$ ?

In this paper, by using a new mountain pass theorem, the perturbation of Newtonian potential $V(q)$, some estimates about the Lagrangian
action, and some other known results, we will give a positive answer to the above question. More precisely, without the winding number condition, we will prove that for planar 3-body problem, beside the variational
minimizer of the Lagrangian action restricted on the loop space $\overline{\Lambda_1}$, there exists another trajectory $q$, which also just forms an Euler's central configuration at any instant.

\vskip2mm\par
The organization of this paper is as follows: In Section 2, we give the main result. Section 3 is devoted to introducing some useful lemmas,  and in Section 4, we prove the main result.

\section{Main result}\label{sec2}\setcounter{equation}{0}
In what follows, for any given positive masses $m_1$, $m_2$ and $m_3$, the configuration space of 3-body problem in $\mathbb{R}^{2}$ is described as
\begin{eqnarray*}
G=\{(x_1,\,x_2,\,x_3)\in(\mathbb{R}^{2})^{3},\,\sum_{i=1}^{3}m_{i}x_{i}=0,\,and \,\,x_i\neq x_j \,\,for \,\,i\neq j\}.
\end{eqnarray*}

For the periodic $T>0$ and the special parameter $\lambda_0$ in (\ref{1.2})-(\ref{1.3}), we also define a new loop space without the winding number condition as the following
\begin{eqnarray}\label{2.1}
\Lambda_2&=&\{q=(q_1,\,q_2,\,q_3)\in G, \,q_i\in W^{1,2}(\mathbb{R}/T\mathbb{Z},\,\mathbb{R}^{2}),
\nonumber \\  && q_{3}(t)-q_{1}(t)=\lambda_{0} (q_{2}(t)-q_{1}(t)), q(t+\frac{T}{2})=-q(t)\}.
\end{eqnarray}
In this paper, we choose $T$ such that
\begin{eqnarray}\label{2.2}
0<T<\Big(\frac{a(\lambda_0)}{2b(\lambda_0)}\Big)^{\frac12}\pi,
\end{eqnarray}
where
\begin{eqnarray}\label{2.3}
\cases{a(\lambda_0)=\frac{1}{\sum_{i=1}^{3}m_{i}}\[m_1m_2+m_1m_3\lambda^{2}_{0}+m_2m_3(1-\lambda_0)^{2}\],
\cr
b(\lambda_0)=m_1m_2+m_1m_3\lambda^{-1}_{0}+m_2m_3(1-\lambda_{0})^{-1}.}
\end{eqnarray}
From \cite[Pages 159-160]{Zhang2018}, we know:
 \begin{itemize}
  \item[$(i)$] For any given positive masses $m_1, \,m_2,\,m_3$, there exists a unique $0<\lambda_{0}<1$ satisfying equation (\ref{1.3});
 \item[$(ii)$] For the $\lambda_0$, any instant $t$ and the non-collision solution $q=(q_1,\,q_2,\,q_3)$, if the solution $q$ satisfies that collinear condition, i.e. $q_{3}(t)-q_{1}(t)=\lambda_{0} (q_{2}(t)-q_{1}(t))$ holds, then $q$ is just the Euler's central configuration at any instant.
\end{itemize}

\begin{remark}\label{rem2.1}
From system (\ref{1.1}), it is easy to see that the non-collision condition $q_i\neq q_j \,for \,i\neq j$ in our new loop space $\Lambda_2$ is natural.
\end{remark}

\begin{remark}\label{rem2.2}
Obviously, combining the above $(ii)$, the periodic and non-collision solution of system (\ref{1.1}) in $\Lambda_2$, is Euler's central configuration at any instant.
\end{remark}

\begin{remark}\label{rem2.3}
From system (\ref{1.1}), we define the Lagrangian action on $\Lambda_2$:
\begin{eqnarray}\label{2.4}
f(q)=\int_{0}^{T}\Big[\frac{1}{2}\sum_{i=1}^{3}m_{i}|\dot{q}_{i}|^{2}+\sum_{1\leq i< j\leq3}\frac{m_{i}m_{j}}{|q_{i}-q_{j}|}\Big]dt.
\end{eqnarray}
In order to get the existence of periodic solution for system (\ref{1.1}), a well known technique is to find the critical point of $f(q)$. Sine the new loop space $\Lambda_2$ is not complete and has special collinear geometry structure with the ratio $\lambda_0$,
there is no obvious extension of the critical theory directly to this case. Thus new ideas are
needed to approach non-complete space with this special collinear geometry structure. In this paper,  we will use the perturbation of the potential $V(q)=-\sum_{1\leq i<j\leq 3}\frac{m_j m_i}{|q_i-q_j|}$ to establish the new existence of non-collision solution for system (\ref{1.1}) such that the three bodies form an Euler's central configuration at any instant.
\end{remark}

We are now in a position to state our main result.
\begin{theorem} \label{the2.1}
For planar 3-body problem (\ref{1.1}) with any given positive masses $m_1, \,m_2,\,m_3$, there exists a trajectory $\tilde{q}$ in the loop space $\Lambda_2$, which is different from the minimizer of $f(q)$ on $\overline{\Lambda_1}$, is also an Euler's central configuration at any instant.
\end{theorem}

\begin{remark}\label{rem2.4}
(i) Since $\tilde{q}\in \Lambda_2$, we  do not need the winding number condition: $deg(q_i-q_j)\neq0 \,(i\neq j)$.

(ii) For planar 3-body problem, by using a direct variational method, Zhang and Zhou \cite{Zhang2004} proved that the minimizer of $f(q)$ on $\overline{\Lambda_1}$, is just an Euler's central configuration at any instant. But now for planar 3-body problem, by using
a different philosophy, i.e. a new mountain pass theorem, the perturbation of Newtonian potential $V(q)$, some estimates about the Lagrangian
action and some other known results, we proved that there exists another trajectory $\tilde{q}$ such that the three bodies form an Euler's central configuration at any instant, i.e. the minimizer of $f(q)$ on the suitable loop space $\overline{\Lambda_1}$, is not the unique solution.
\end{remark}

\section{Useful lemmas}\setcounter{equation}{0}\setcounter{equation}{0}
In 2019, Ding, Wei and Zhang \cite{Ding2019} obtained the following extension of the mountain pass theorem.

\begin{lemma} \cite[Theorem 1.4]{Ding2019} \label{lem3.1} Let $X$ be a Hilbert space,
$\tilde{f}\in C^{2}(X, \mathbb{R})$, $q^{(e)}, q^{(e_1)}\in X$ and
 $r>0$ such that $0<\|q^{(e_1)}\|_X<r$ and $\|q^{(e)}\|_X>r$, and
$\tilde{f}(\theta)<\tilde{f}(q^{(e)})=\tilde{f}(q^{(e_1)})$.
Then, for each small enough $\varepsilon>0$, there exists $\hat{q}\in X$ such that
\begin{itemize}
\item[$(i)$] $\hat{c}-2\varepsilon\leq \tilde{f}(\hat{q})\leq \hat{c}+2\varepsilon$;
\item[$(ii)$] $\|\tilde{f}^{\prime}(\hat{q})\|_X<2\varepsilon$,
\end{itemize}
where $\hat{c}:=\inf_{\gamma\in\hat{\Gamma}}\max_{t\in[0, 1]}\tilde{f}\big(\gamma(t)\big)$ and
\begin{eqnarray*}
\hat{\Gamma}:=\{\gamma\in C\big([0, 1], X\big):\gamma(0)=\theta, \gamma(\frac{1}{2})=q^{(e_1)},\gamma(1)=q^{(e)}\}.
\end{eqnarray*}
\end{lemma}

\begin{remark}\label{rem3.1}
Let $c_0:=\inf_{\|q\|_X=r}\tilde{f}(q)$ and let $c_1:=\max\{\tilde{f}(q^{(\theta)}), \,\,\tilde{f}(q^{(e)})\}$. For the well-known mountain pass type theorems in references \cite{Ambrosetti1973,Br1980,Chang1993,Pucci1985,Qi1987,Willem1996}, we see that $c_0\geq c_1$.
But in Lemma \ref{lem3.1}, it is independent of $c_0$, which implies that Lemma \ref{lem3.1} holds not only for $c_0\geq c_1$, but also holds for $c_0<c_1$.
\end{remark}

Let $\Lambda_2$ be given by (\ref{2.1}) and let $\overline{\Lambda_2}$ be the closure of $\Lambda_2$ (with the same norm of $\Lambda_2$). For $q\in\Lambda_2$, since $q_{i}(t+\frac{T}{2})=-q_{i}(t)\,(i=1,2,3)$,
then $\int_{0}^{T}q_i(t)dt=0\,\, (i=1,2,3)$. With the help of the Poincar\'{e} inequality, for $q\in\Lambda_2$ the norm in $\Lambda_2$ can be defined by
$$
\|q\|=\[\int_{0}^{T}\(\sum_{i=1}^{3}m_{i}|\dot{q}_{i}|^{2}\)dt\]^{1/2}.
$$

Employing Lemma \ref{lem3.1}, we have the following result which also holds without the restriction of $c_0\geq c_1$.
\begin{lemma}\label{lem3.2} Let $\tilde{f}\in C^{2}(\Lambda_2, \mathbb{R})$, $q^{(\theta)}, \,q^{(e)}, q^{(e_1)}\in \Lambda_2$
and $r>0$ such that $0<\|q^{(\theta)}\|<\|q^{(e_1)}\|<r$ and $\|q^{(e)}\|>r$, and
$\tilde{f}(q^{(\theta)})<\tilde{f}(q^{(e)})=\tilde{f}(q^{(e_1)})$.
Then, for each small enough $\varepsilon>0$, there exists $\hat{q}\in \Lambda_2$ such that
\begin{itemize}
\item[$(i)$] $\hat{c}-2\varepsilon\leq\tilde{f}(\hat{q})\leq \hat{c}+2\varepsilon$;
\item[$(ii)$] $\|\tilde{f}^{\prime}(\hat{q})\|<2\varepsilon$,
\end{itemize}
where $\hat{c}:=\inf_{\gamma\in\hat{\Gamma}}\max_{t\in[0, 1]}\tilde{f}\big(\gamma(t)\big)$ and
\begin{eqnarray*}
\hat{\Gamma}:=\{\gamma\in C\big([0, 1], \Lambda_2\big):\gamma(0)=q^{(\theta)}, \gamma(\frac{1}{2})=q^{(e_1)},\gamma(1)=q^{(e)}\}.
\end{eqnarray*}
\end{lemma}
\textbf{Proof.} The proof can be followed by Lemma \ref{lem3.1}. To avoid unnecessary repetition, we only give the frame of the proof.

\textbf{Step 1.} From the calculations in quantitative deformation lemma \cite[Lemma 2.1]{Ding2019} (also see \cite{Ambrosetti1973,Ambrosetti1992}), it is clearly to find that if we change the Hilbert space $X$ into the above loop space $\Lambda_2$, the quantitative deformation lemma holds as well, i.e.,

Let $\Lambda_2$ be defined as (\ref{2.1}), and $\varepsilon$ be a small enough positive number. Let
$\tilde{f}\in C^{2}(\Lambda_2, \mathbb{R})$, $s\in \mathbb{R}$. Assume that
\begin{eqnarray*}
\|\tilde{f}^{\prime}(q)\|\geq2\varepsilon, \quad \forall \  q\in \tilde{f}^{-1}([s-2\varepsilon,
s+2\varepsilon]).
\end{eqnarray*}
Then there exists $\eta\in$ $C(\Lambda_2, \Lambda_2)$, such that
\begin{itemize}
\item[$(i)$] $\eta(q)=q$, $\forall \ q\notin \tilde{f}^{-1}\big([s-2\varepsilon,
s+2\varepsilon]\big)\backslash D$, where $D$ is any
subset of $\Lambda_2$ satisfying
$D\subset \tilde{f}^{-1}\big([s-\frac{1}{3}\varepsilon^{2},
s+\frac{1}{3}\varepsilon^{2}]\big)$;
\item[$(ii)$] $\eta\big(\tilde{f}^{-1}[s+\frac{1}{2}\varepsilon^{2}, s+\varepsilon^{2}]\big)\subset \tilde{f}^{-1}\big([s-\frac{3}{2}\varepsilon^{2},
s-\frac{1}{2}\varepsilon^{2}]\big)$.
\end{itemize}

\textbf{Step 2.} By \textbf{Step 1},  the quantitative deformation lemma \cite[Lemma 2.1]{Ding2019} is true for $\Lambda_2$, so  all proof manipulations appealing to $X$ in Lemma \ref{lem3.1} is adapted to $\Lambda_2$ here again. Therefore, we finish the proof. $\Box$

\begin{remark}\label{rem3.2}
In Lemma \ref{lem3.2}, if $\tilde{f}(q)$ satisfies the $(P.S.)$ condition in $\Lambda_2$, then there exists $\hat{q}\in \overline{\Lambda_2}$ such that $\tilde{f}(\hat{q})= \hat{c}$,
 which implies $\hat{q}$ is the critical point of $\tilde{f}$ in $\overline{\Lambda_2}$.
\end{remark}

Since planar 3-body system (\ref{1.1}) is equivalent to the following system
\begin{eqnarray*}
m_i\ddot{q}_i+\nabla_{q_{i}} V(q)=0, \quad q_{i}\in\mathbb{R}^{2}, \quad \ i=1,2,3,
\end{eqnarray*}
where
\begin{eqnarray*}
V(q)=-\sum_{1\leq i< j\leq3}\frac{m_im_j}{|q_i-q_j|}, \ \ \ \ q=(q_{1},\,q_{2},\,q_{3}),
\end{eqnarray*}
then in order to study the new existence of non-collision solution for planar 3-body system (\ref{1.1}) that the three bodies form an Euler's central configuration at any instant, firstly we consider the following system, i.e. the perturbation of fixed energy system:
\begin{eqnarray}\label{3.1}
\cases{m_i\ddot{q}_i+\nabla_{q_{i}} V_{\tilde{\varepsilon}}(q)=0 , \quad q_{i}\in\mathbb{R}^{2}, \quad i=1,2,3,\cr
\frac{1}{2}\sum_{1\leq i\leq 3} m_i|\dot{q}_i(t)|^{2}+ V_{\tilde{\varepsilon}}(q)=h,}
\end{eqnarray}
where $\tilde{\varepsilon}>0$ is a constant, $h$ represents the fixed energy and
\begin{eqnarray}\label{3.2}
V_{\tilde{\varepsilon}}(q)=-\sum_{1\leq i< j\leq3}\frac{m_{i}m_{j}}{|q_{i}-q_{j}|}+\tilde{\varepsilon}\cdot\frac{1}{h}\cdot\sum_{1\leq i< j\leq3}\frac{m_{i}m_{j}}{|q_{i}-q_{j}|^{2}}.
\end{eqnarray}
From $b(\lambda_0)>0$, we can set $h$ such that
$-b(\lambda_0)/2<h<0$. By $\|q\|=\[\int_{0}^{T}\(\sum_{i=1}^{3}m_{i}|\dot{q}_{i}|^{2}\)dt\]^{1/2}$ with $q_{i}(t)\in\mathbb{R}^{2}$, we consider a new functional
\begin{eqnarray}\label{3.3}
\varphi_{\tilde{\varepsilon}}(q)&=&\frac{1}{2}\|q\|^{2}\int_{0}^{T}\(h-V_{\tilde{\varepsilon}}(q)\)dt\nonumber\\
&=&\frac{1}{2}\|q\|^{2}\int_{0}^{T}\(h+\sum_{1\leq i< j\leq3}\frac{m_{i}m_{j}}{|q_{i}-q_{j}|}- \tilde{\varepsilon}\cdot\frac{1}{h}\cdot\sum_{1\leq i< j\leq3}\frac{m_{i}m_{j}}{|q_{i}-q_{j}|^{2}}\)dt.
\end{eqnarray}
Observing that $\Lambda_2\subseteq \{q(t)=(q_1,\,q_2,\,q_3), \,q_i\in W^{1,2}(\mathbb{R}/T\mathbb{Z},\,\mathbb{R}^{2}),
\,q(t+\frac{T}{2})=-q(t),\, q_i\neq q_j \,\,for \,\,i\neq j\}\subseteq W^{1,2}(\mathbb{R}/T\mathbb{Z},\,\mathbb{R}^{2})$, then we have the following lemmas.

\begin{lemma}\label{lem3.3}\cite[Lemma 1]{Ambrosetti1992}
For any $\tilde{\varepsilon}>0$, let $q_{\tilde{\varepsilon}}\in \Lambda_2$ be such that $\varphi_{\tilde{\varepsilon}}^{\prime}(q_{\tilde{\varepsilon}})=0$ and $\|q_{\tilde{\varepsilon}}\|>0$ and set
\begin{eqnarray*}
\omega^{2}_{\tilde{\varepsilon}}=\frac{\int_{0}^{T}\nabla V_{\tilde{\varepsilon}}(q_{\tilde{\varepsilon}})q_{\tilde{\varepsilon}}}{\|q_{\tilde{\varepsilon}}\|^{2}}>0.
\end{eqnarray*}
Then for any fixed energy $-b(\lambda_0)/2<h<0$, $\tilde{q}_{\tilde{\varepsilon}}(t):=q_{\tilde{\varepsilon}}(\omega_{\tilde{\varepsilon}}t)$ is a periodic and non-collision solution of system (\ref{3.1}).
\end{lemma}

\begin{lemma}\label{lem3.4}
If $q^{(n)}\rightharpoonup q$ in $\Lambda_2$ and $q\in \partial \Lambda_2$, then $\int_{0}^{T}V_{\tilde{\varepsilon}}(q^{(n)})dt\rightarrow-\infty$.
\end{lemma}
\textbf{Proof.} Noting that $q^{(n)}\rightharpoonup q\in\partial \Lambda_2\subset W^{1,2}(\mathbb{R}/T\mathbb{Z},\,\mathbb{R}^{2})$, by Sobolev's compact embedding theorem (or see \cite[Proposition 1.2]{Mawhin1989}), we know $q^{(n)}\rightarrow q$
uniformly on $[0,T]$. Since $q\in\partial \Lambda_2$, there exists $t^{\ast}\in[0,T]$, $q_{l}$ and $q_{r}$ ($l,r\in\{1,2,3\}$) such
that $q_{l}(t^{\ast})-q_{r}(t^{\ast})=0$.  Therefore, there exists $\delta>0$ such that $|q^{(n)}_{l}(t)-q^{(n)}_{r}(t)|<\tilde{\varepsilon}$ for all
$t\in[t^{\ast}-\delta,\,t^{\ast}+\delta]$ if $n$ is large enough.

By (\ref{3.2}), we have
\begin{eqnarray}\label{3.4}
\int_{0}^{T}V_{\tilde{\varepsilon}}(q^{(n)})dt\leq\int_{0}^{T}\frac{\tilde{\varepsilon}}{h}\sum_{1\leq i< j\leq3}\frac{m_{i}m_{j}}{|q^{(n)}_{i}-q^{(n)}_{j}|^{2}}dt\leq\int_{t^{\ast}-\delta}^{t^{\ast}+\delta}
\frac{\tilde{\varepsilon}}{h}\frac{m_{l}m_{r}}{|q^{(n)}_{l}-q^{(n)}_{r}|^{2}}dt
\end{eqnarray}
for $h<0$.

If $q_{l}-q_{r}\equiv0$ for all $t\in[t^{\ast}-\delta,\,t^{\ast}+\delta]$, then we have
\begin{eqnarray*}
\int_{t^{\ast}-\delta}^{t^{\ast}+\delta}\frac{m_{l}m_{r}}{|q^{(n)}_{l}(t)-q^{(n)}_{r}(t)|^{2}}dt\rightarrow+\infty,
\end{eqnarray*}
because $q^{(n)}_{l}-q^{(n)}_{r}\rightarrow 0$
uniformly on $[t^{\ast}-\delta,\,t^{\ast}+\delta]$. So we can
assume that $q_{l}(t^{\ast}+\delta)-q_{r}(t^{\ast}+\delta)\neq 0$ and we find
\begin{eqnarray}\label{3.5}
\[\log|q_{l}^{(n)}(t)-q_{r}^{(n)}(t)|\]_{t^{\ast}}^{t^{\ast}+\delta}&=&
\int_{t^{\ast}}^{t^{\ast}+\delta}\frac{|\dot{q}^{(n)}_{l}-\dot{q}^{(n)}_{r}|}{|q^{(n)}_{l}-q^{(n)}_{r}|}dt
\nonumber\\&\leq&
\|\dot{q}^{(n)}_{l}-\dot{q}^{(n)}_{2}\|_{L^2}\[\int_{t^{\ast}}^{t^{\ast}+\delta}
\frac{1}{|q^{(n)}_{l}-q^{(n)}_{r}|^2}dt\]^{\frac{1}{2}}.
\end{eqnarray}
Noting that $q^{(n)}_{l}-q^{(n)}_{r}$ converges weakly to $q_{l}-q_{r}$, then $\|q^{(n)}_{l}-q^{(n)}_{r}\|_{L^{2}}+\|\dot{q}^{(n)}_{l}-\dot{q}^{(n)}_{r}\|_{L^{2}}$ is bounded.
So $\|\dot{q}^{(n)}_{l}-\dot{q}^{(n)}_{r}\|_{L^{2}}$ is bounded. Since $q^{(n)}_{l}(t^{\ast})-q^{(n)}_{r}(t^{\ast})\rightarrow0$ while $q^{(n)}_{l}(t^{\ast}+\delta)-q^{(n)}_{r}(t^{\ast}+\delta)\rightarrow q_{l}(t^{\ast}+\delta)-q_{r}(t^{\ast}+\delta)\neq0$,
then (\ref{3.5}) and $\[\log|q_{l}^{(n)}(t)-q_{r}^{(n)}(t)|\]_{t^{\ast}}^{t^{\ast}+\delta}\rightarrow+\infty$ immediately imply that
\begin{eqnarray}\label{3.6}
+\infty\leq\liminf_{n\rightarrow\infty} \int_{t^{\ast}}^{t^{\ast}+\delta}\frac{1}{|q^{(n)}_{l}-q^{(n)}_{r}|^{2}}dt.
\end{eqnarray}
Then, by (\ref{3.4}), (\ref{3.6}) and Fatou's lemma, it yields that
\begin{eqnarray*}
\limsup_{n\rightarrow\infty}\int_{0}^{T}V_{\tilde{\varepsilon}}(q^{(n)})dt&\leq&
\limsup_{n\rightarrow\infty} \[\int_{t^{\ast}-\delta}^{t^{\ast}+\delta}
\frac{\tilde{\varepsilon}}{h}\frac{m_{l}m_{r}}{|q^{(n)}_{l}-q^{(n)}_{r}|^{2}}dt\]
\nonumber\\&=&
\frac{\tilde{\varepsilon}}{h}\liminf_{n\rightarrow\infty} \[\int_{t^{\ast}-\delta}^{t^{\ast}+\delta}
\frac{m_{l}m_{r}}{|q^{(n)}_{l}-q^{(n)}_{r}|^{2}}dt\]
=-\infty.
\end{eqnarray*}
From this, we complete the proof. $\Box$

\begin{lemma}\cite[Lemma 2, conclusion (i)]{Ambrosetti1992}\label{lem3.5}
There exist $\rho,\,\beta>0$ such that $\varphi_{\tilde{\varepsilon}}(q)\geq\beta$ for all $\tilde{\varepsilon}>0$ and all $q\in \Lambda_2$, $\|q\|=\rho$.
\end{lemma}

From the definition of the functional $\varphi_{\tilde{\varepsilon}}(q)$ and the loop space $\Lambda_2$, we know that $\varphi_{\tilde{\varepsilon}}\in C^{2}(\Lambda_2, \mathbb{R})$. Moreover, we can prove the following lemma.
\begin{lemma}\label{lem3.6}
There exist $q^{(\theta)}, \,q^{(e)}, \,q^{(e_1)}\in \Lambda_2$,
 $r>0$ such that $0<\|q^{(\theta)}\|<\|q^{(e_1)}\|<r$ and $\|q^{(e)}\|>r$, and
$\varphi_{\tilde{\varepsilon}}(q^{(\theta)})<\varphi_{\tilde{\varepsilon}}(q^{(e)})=\varphi_{\tilde{\varepsilon}}(q^{(e_1)})$.
\end{lemma}
\textbf{Proof.} We divide the proof into two steps.

\textbf{Step 1.} We search for $q^{(e)}, \,q^{(e_1)}\in \Lambda_2$ such that $\|q^{(e_1)}\|<\|q^{(e)}\|$ and $\varphi_{\tilde{\varepsilon}}(q^{(e)})=\varphi_{\tilde{\varepsilon}}(q^{(e_1)})$.

 From (\ref{3.3}), $\|q\|=\[\int_{0}^{T}\(\sum_{i=1}^{3}m_{i}|\dot{q}_{i}|^{2}\)dt\]^{\frac{1}{2}}$ with $q_{i}(t)\in\mathbb{R}^{3}$, and $q_{3}(t)-q_{1}(t)=\lambda_{0}(q_{2}(t)-q_{1}(t))$, we have
\begin{eqnarray}\label{3.7}
\varphi_{\tilde{\varepsilon}}(q)&=&\frac{1}{2}\|q\|^{2}\int_{0}^{T}\Big[h+\sum_{1\leq i< j\leq3}\frac{m_{i}m_{j}}{|q_{i}-q_{j}|}
\nonumber\\ &&+\tilde{\varepsilon} \cdot(-\frac{1}{h})\cdot\sum_{1\leq i< j\leq3}\frac{m_{i}m_{j}}{|q_{i}-q_{j}|^{2}}\Big]dt\nonumber\\
&=&\frac{1}{2}\int_{0}^{T}\(\sum_{i=1}^{3}m_{i}|\dot{q}_{i}|^{2}\)dt\int_{0}^{T}\Big[h+\sum_{1\leq i< j\leq3}\frac{m_{i}m_{j}}{|q_{i}-q_{j}|}\nonumber\\ &&+\tilde{\varepsilon} \cdot(-\frac{1}{h})\cdot\sum_{1\leq i< j\leq3}\frac{m_{i}m_{j}}{|q_{i}-q_{j}|^{2}}\Big]dt\nonumber\\
&=&\frac{1}{2}\int_{0}^{T}\(\sum_{i=1}^{3}m_{i}|\dot{q}_{i}|^{2}\)dt
\int_{0}^{T}\[h+(m_1m_2+m_1m_3\lambda^{-1}_{0}+\frac{m_2m_3}{1-\lambda_{0}}\frac{1}{|q_{2}-q_{1}|})
\nonumber\\
&&+\tilde{\varepsilon} \cdot(-\frac{1}{h})\cdot(m_1m_2+\frac{m_1m_3}{\lambda^2_0}+\frac{m_2m_3}{(1-\lambda_{0})^2})\frac{1}{|q_{2}-q_{1}|^{2}}\]dt.
\end{eqnarray}
Let $q^{(e)}_{i}(t)\in\mathbb{R}^{2}$ with $i=1,2$ such that $|q^{(e)}_2(t)-q^{(e)}_{1}(t)|\equiv1$. We set $q^{(e)}=(q^{(e)}_1,q^{(e)}_2,0)$, $ q^{(e_1)}=(q^{(e_1)}_1,q^{(e_1)}_2,0)=\mu\cdot q^{(e)}$,
 and $\mu=[-b(\lambda_0)-h]/h$. Observing that $-b(\lambda_0)/2<h<0$, one computes that $\mu>1$. From $b(\lambda_0)=m_1m_2+m_1m_3\lambda^{-1}_{0}+m_2m_3(1-\lambda_{0})^{-1}$, we have
\begin{eqnarray*}
\varphi_{\tilde{\varepsilon}}(q^{(e)})=\varphi_{\tilde{\varepsilon}}(q^{(e_1)})
\Longleftrightarrow \int_{0}^{T}\Big[h\mu^{2}+b(\lambda_0)\mu-\Big(b(\lambda_0)+h\Big)\Big]dt=0.
\end{eqnarray*}
 Since $\mu=[-b(\lambda_0)-h]/h$ and $-b(\lambda_0)/2<h<0$, we get
\begin{eqnarray*}
 \int_{0}^{T}\Big[h\mu^{2}+b(\lambda_0)\mu-\Big(b(\lambda_0)+h\Big)\Big]dt=0,
\end{eqnarray*}
 which implies that there exists $\|q^{(e_1)}\|<\mu\|q^{(e_1)}\|=\|q^{(e)}\|$ such that $\varphi_{\tilde{\varepsilon}}(q^{(e)})$=$\varphi_{\tilde{\varepsilon}}(q^{(e_1)})$.

\textbf{Step 2.} We search for the $q^{(\theta)}$ such that $\varphi_{\tilde{\varepsilon}}(q^{(\theta)})<\varphi_{\tilde{\varepsilon}}(q^{(e)})=\varphi_{\tilde{\varepsilon}}(q^{(e_1)})$.

Denote $q^{(E)}=x\cdot q^{(e)}=x\cdot(q^{(e)}_1,q^{(e)}_2,0)$ with $x>0$, and $|q^{(e)}_2(t)-q^{(e)}_{1}(t)|\equiv1$,
then by (\ref{3.7}), we have
\begin{eqnarray}\label{3.8}
\varphi_{\tilde{\varepsilon}}(q^{(E)})&=&\frac{1}{2}\int_{0}^{T}\(\sum_{i=1}^{3}m_{i}|\dot{q}^{(E)}_{i}|^{2}\)dt
\nonumber\\ &&\times\int_{0}^{T}\[h+\frac{s}{|q^{(E)}_{2}-q^{(E)}_{1}|}
+\tilde{\varepsilon} \cdot(-\frac{1}{h})\cdot s\cdot\frac{1}{|q^{(E)}_{2}-q^{(E)}_{1}|^{2}}\]dt\nonumber\\
&=&\frac{x^{2}}{2}\int_{0}^{T}\(\sum_{i=1}^{3}m_{i}|\dot{q}^{(e)}_{i}|^{2}\)dt\nonumber\\ &&\times\int_{0}^{T}\[h+\frac{s}{|x|\cdot|q^{(e)}_{2}-q^{(e)}_{1}|}
-\frac{\tilde{\varepsilon}}{h}\cdot\frac{s}{x^{2}\cdot|q^{(e)}_{2}-q^{(e)}_{1}|^{2}}\]dt.
\end{eqnarray}
Thus employing (\ref{3.8}) and $|q^{(e)}_2(t)-q^{(e)}_{1}(t)|\equiv1$, we can define the function $g(x)$ as follows:
\begin{eqnarray*}
g(x)&\triangleq&\varphi_{\tilde{\varepsilon}}(q^{(E)})=\frac{x^{2}}{2}\int_{0}^{T}\(\sum_{i=1}^{3}m_{i}|\dot{q}^{(e)}_{i}|^{2}\)dt\int_{0}^{T}\[h+\frac{s}{x}
+\tilde{\varepsilon} \cdot(-\frac{1}{h})\cdot\frac{s}{x^{2}}\]dt\nonumber\\
&=&\frac{Th\int_{0}^{T}(\sum_{i=1}^{3}m_{i}|\dot{q}^{(e)}_{i}|^{2})dt}{2} x^{2}+ \frac{sT\int_{0}^{T}(\sum_{i=1}^{3}m_{i}|\dot{q}^{(e)}_{i}|^{2})dt}{2}x
\nonumber\\
&&-\frac{\tilde{\varepsilon} sT\int_{0}^{T}(\sum_{i=1}^{3}m_{i}|\dot{q}^{(e)}_{i}|^{2})dt}{2h}.
\end{eqnarray*}
Note that $\frac{Th\int_{0}^{T}(\sum_{i=1}^{3}m_{i}|\dot{q}^{(e)}_{i}|^{2})dt}{2}<0$, $x>0$ and
\begin{eqnarray*}
\cases{
q^{(E)}=q^{(e)}=x\cdot q^{(e)},\,\, where \,\, x=1, \cr
q^{(E)}=q^{(e_1)}=x\cdot q^{(e)},\,\,where \,\,x=\mu>1,\cr
\varphi_{\tilde{\varepsilon}}(q^{(e)})=g(1)=\varphi_{\tilde{\varepsilon}}(q^{(e_1)})=g(\mu).}
\end{eqnarray*}
Then we conclude that there exists $0<x<<1$ such that
\begin{eqnarray*}
\cases{q^{(\theta)}=x\cdot q^{(e)},\,\, where \,\, 0<x<<1, \cr
\|q^{(\theta)}\|<\|q^{(e_1)}\|<\|q^{(e)}\|,\cr
\varphi_{\tilde{\varepsilon}}(q^{(\theta)})<\varphi_{\tilde{\varepsilon}}(q^{(e)})=\varphi_{\tilde{\varepsilon}}(q^{(e_1)}).}
\end{eqnarray*}
By now, we complete the proof of Lemma \ref{lem3.6}. $\Box$

In \cite{Ambrosetti1992}, the authors assumed that $\Omega=\mathbb{R}^{2}\backslash\{0\}$, and the potential $V(q)=\sum_{1\leq i<j\leq N}V_{ji}(q_j-q_i)$ with $V_{ji}$ satisfying
\begin{itemize}
\item [$(V_1)$] $V_{ji}(\xi)=V_{ij}(\xi)$, \,\,$\forall\,\,\xi\in\Omega$;
\item [$(V_2)$] $\exists\,\alpha\in[1,2)$ such that $\nabla\,V_{ji}(\xi)\xi\geq-\alpha V_{ji}(\xi)>0$, $\forall \,\xi\in\Omega$;
\item [$(V_3)$] $\exists\,\delta^{\prime}\in(0,\,2)$ and $r>0$ such that $\nabla\,V_{ji}(\xi)\xi\leq-\delta^{\prime}V_{ji}(\xi)$ for all $0<|\xi|\leq r$;
\item [$(V_4)$] $V_{ji}(\xi)\rightarrow0$ as $|\xi|\rightarrow+\infty$.
\end{itemize}

Clearly, the potential $V(q)=-\sum_{1\leq i<j\leq 3}\frac{m_j m_i}{|q_j-q_i|}$ in this paper satisfies conditions $(V_1)-(V_4)$. We set
\begin{eqnarray*}
\Lambda_0=\{q(t)=(q_1,\,q_2,\,q_3), \,q_i\in W^{1,2}(\mathbb{R}/T\mathbb{Z},\,\mathbb{R}^{2}),
\,q(t+\frac{T}{2})=-q(t),\,
q_i\neq q_j \,\,for \,\,i\neq j\}.
\end{eqnarray*}
With the aid of Lemmas \ref{lem3.2} and \ref{lem3.6}, there exists a sequence $\{q^{(n)}\}\subseteq\Lambda_2\subseteq\Lambda_0$, such that $\varphi_{\tilde{\varepsilon}}(q^{(n)})\rightarrow \hat{c}$,
and $\varphi^{\prime}_{\tilde{\varepsilon}}(q^{(n)})\rightarrow0$. Combining Lemma \ref{lem3.5},
there exist constants $M_1$ and $\beta$ such that $0<\beta\leq\varphi_{\tilde{\varepsilon}}(q^{(n)})\leq M_1$ and $\varphi^{\prime}_{\tilde{\varepsilon}}(q^{(n)})\rightarrow0$. Then employing
$h<0$ and Lemma \ref{lem3.4}, all the conditions of \cite[Lemma 5]{Ambrosetti1992} are satisfied. Thus we have

\begin{lemma}\label{lem3.7}\cite[Lemma 5]{Ambrosetti1992}
If $q^{(n)}\in\Lambda_2$ is such that $0<\varphi_{\tilde{\varepsilon}}(q^{(\theta)})<\varphi_{\tilde{\varepsilon}}(q^{(n)})\leq M_1$ where $M_1$ is a positive constant, and
$\varphi_{\tilde{\varepsilon}}^{\prime}(q^{(n)})\rightarrow0$, then (up to a subsequence) $q^{(n)}\rightarrow q^{\ast}\in \Lambda_2$.
\end{lemma}

In Lemma \ref{lem3.2}, take $\tilde{f}(q)=\varphi_{\tilde{\varepsilon}}(q)$. Then we have
\begin{lemma}\label{lem3.8}
There exists $\varepsilon_{0}>0$ such that for any $\tilde{\varepsilon}\in(0,\,\varepsilon_{0})$, there is $q_{\tilde{\varepsilon}}$ satisfied that $q_{\tilde{\varepsilon}}$ is the critical point of $\varphi_{\tilde{\varepsilon}}$ in $\Lambda_{2}$. Moreover, there exist $a,\,b>0$
such that $0<a\leq \|q_{\tilde{\varepsilon}}\|\leq b$ holds for any $\tilde{\varepsilon}\in(0,\,\varepsilon_{0})$.
\end{lemma}
\textbf{Proof.} By Lemmas \ref{lem3.6} and \ref{lem3.7}, we know all the assumptions of Lemma \ref{lem3.2} are satisfied. Then by Lemma \ref{lem3.2} and Remark \ref{rem3.2},
there exists a critical point $q_{\tilde{\varepsilon}}\in\Lambda_{2}$ of $\varphi_{\tilde{\varepsilon}}$. The rest of the proof is the same as the proof of Lemma 6 in \cite{Ambrosetti1992}, so we omit the details.

\begin{remark}\label{rem3.3}
(i) The same conclusion of Lemma \ref{lem3.8} was obtained in \cite{Ambrosetti1992} with the condition of $q_{\tilde{\varepsilon}}\in \Lambda_0$,
\begin{eqnarray*}
\Lambda_0=\{q(t)=(q_1,\,q_2,\,q_3), \,q_i\in W^{1,2}(\mathbb{R}/T\mathbb{Z},\,\mathbb{R}^{2}),
\,q(t+\frac{T}{2})=-q(t),\, q_i\neq q_j \,\,for \,\,i\neq j\},
\end{eqnarray*}
but in Lemma \ref{lem3.8}, $q_{\tilde{\varepsilon}}\in \Lambda_{2}$.

(ii) In 1992, Ambrosetti and Zelati used the conclusion $(ii)$ of Lemma 2 in \cite{Ambrosetti1992} (i.e. there exist $\varepsilon_{0}>0$, $q^{(\theta)},\,q^{(e)}\in\Lambda_{2}$
with $\|q^{(\theta)}\|<\rho<\|q^{(e)}\|$, and positive constant $\rho,\,\beta$ which satisfies Lemma \ref{lem3.5}, such that $\varphi_{\tilde{\varepsilon}}(q^{(\theta)})<\beta,\,\varphi_{\tilde{\varepsilon}}(q^{(e)})<\beta$ for any $\tilde{\varepsilon}\in(0,\,\varepsilon_{0})$), to obtain the same conclusion of Lemma \ref{lem3.8} in $\Lambda_{2}$. But now since the loop space $\Lambda_{2}$ has geometry structure $q_{3}(t)-q_{1}(t)=\lambda_{0} (q_{2}(t)-q_{1}(t))$ where $\lambda_{0}$ satisfies
\begin{eqnarray*}
\frac{m_{3}\lambda_{0}^{-2}+m_2}{m_{3}\lambda_{0}+m_2}-\frac{m_{3}(1-\lambda_{0})^{-2}+m_1}{m_{3}(1-\lambda_{0})+m_1}=0,
\end{eqnarray*}
then for $q\in\Lambda_{2}$, the method in \cite{Ambrosetti1992} is invalid. More precisely, it is very difficult to obtain the conclusion $(ii)$ of Lemma 2 in \cite{Ambrosetti1992}, which implies that it is very difficult to verify $c_0\geq c_1$ $(c_0:=\inf_{\|q\|=r}\varphi_{\tilde{\varepsilon}}(q)$, $c_1:=\max\{\varphi_{\tilde{\varepsilon}}(q^{(\theta)}), \,\,\tilde{f}(q^{(e)})\})$ in the known mountain pass type theorems, but fortunately, Lemma \ref{lem3.2} holds without the restriction of $c_0\geq c_1$.
\end{remark}

Moreover, we need the following Lemma \ref{lem3.9}, and some estimates about the Lagrangian
action of 2-body problem:
\begin{lemma}\label{lem3.9}\cite[Page 167]{Zhang2018}
$\sum_{1\leq i<j\leq3}m_{i}m_{j}|\dot{q}_{i}-\dot{q}_{j}|^{2}=\sum_{1\leq i\leq3}m_{i}|\dot{q}_{i}|^{2}$.
\end{lemma}

Consider the functional $f(q)$ defined in (\ref{2.4}), then by Lemma \ref{lem3.9}, similar to the method of \cite{Zhang2004}, we have
\begin{eqnarray}\label{3.9}
f(q)&=&\int_{0}^{T}\Big[\frac{1}{2}\sum_{i=1}^{3}m_{i}|\dot{q}_{i}|^{2}+\sum_{1\leq i< j\leq3}\frac{m_{i}m_{j}}{|q_{i}-q_{j}|}\Big]dt\nonumber\\
&=&\int_{0}^{T}\frac{1}{2\sum_{i=1}^{3}m_{i}}\[m_1m_2+m_1m_3\lambda^{2}_{0}+
m_2m_3(1-\lambda_0)^{2}\]|\dot{q}_{1}-\dot{q}_{2}|^{2}\nonumber\\
&&+\[m_1m_2+m_1m_3\lambda^{-1}_{0}+m_2m_3(1-\lambda_{0})^{-1}\]\frac{1}{|q_1-q_2|}dt.
\end{eqnarray}
Let $\lambda_{0}$ be given as (\ref{1.3}), and let $a(\lambda_0)$ and $b(\lambda_0)$ given by (\ref{2.3}). Then, (\ref{3.9}) is equivalent to
\begin{eqnarray}\label{3.10}
f(q)=a(\lambda_0)\int_{0}^{T}\Big[\frac{1}{2}|\dot{q}_{1}-\dot{q}_{2}|^{2}+\frac{b(\lambda_0)/a(\lambda_0)}{|q_1-q_2|}\Big]dt.
\end{eqnarray}
Define the loop space
\begin{eqnarray*}
\Lambda^{\prime}_{1}=\{q=(q_1,\,q_2), \,q_1,q_2\in W^{1,2}(\mathbb{R}/T\mathbb{Z},\,\mathbb{R}^{2}),\,q_1\neq q_2,\,
deg(q_2-q_1)\neq0 \}
\end{eqnarray*}
and
\begin{eqnarray}\label{3.11}
f_{1}(q)=\int_{0}^{T}\Big[\frac{1}{2}|\dot{q}_{1}-\dot{q}_{2}|^{2}+\frac{b(\lambda_0)/a(\lambda_0)}{|q_1-q_2|}\Big]dt.
\end{eqnarray}

\begin{lemma}\label{lem3.10}\cite[Lemma 2.1]{Gordon1977}
The minimizer for $f_1(q)$ on $\overline{\Lambda^{\prime}_{1}}$ is precisely the Keplerian elliptical or collision ejection orbit, and the minimum of the action functional $f_1(q)$ equals to
\begin{eqnarray*}
A=(3\pi)(T/2\pi)^{\frac{1}{3}}\[b(\lambda_0)/a(\lambda_0)\]^{\frac23}=\frac32\cdot(2\pi)^{\frac23}
\[b(\lambda_0)/a(\lambda_0)\]^{\frac23}T^{\frac13}.
\end{eqnarray*}
\end{lemma}

\begin{lemma}\label{lem3.11}\cite[Theorem 1.1]{Zhang2004}
The minimizer of $f(q)$ on $\overline{\Lambda_1}$ is precisely the Euler's central configuration at any instant.
\end{lemma}

\begin{lemma}\label{lem3.12}\cite[Theorem 3.2]{Long2000}
Let $q\in W^{1,2}(\mathbb{R}/T\mathbb{Z},\,\mathbb{R}^{2})$ and $\int^{T}_{0}q(t)dt/T=0$, then
\begin{eqnarray*}
\int^{T}_{0}\[\frac{1}{2}|\dot{q}|^{2}+\frac{C}{|q|}\]dt\geq\frac{3}{2}(2\pi)^{\frac{2}{3}}C^{\frac{2}{3}}T^{\frac{1}{3}}.
\end{eqnarray*}
\end{lemma}

\begin{lemma}\label{lem3.13}\cite[Theorem 8.1.3]{Kuczma2009}
Suppose that $\alpha\leq k(x)\leq\beta$, where $\alpha$ and $\beta$ may be finite or infinite, the range of integration
and the weight function $p(x)$ is finite and positive everywhere, and $\phi^{\prime\prime}(t)$ is positive finite for $\alpha<t<\beta$. Then
\begin{eqnarray*}
\phi(\frac{\int^{\beta}_{\alpha}k(x)p(x)dx}{\int^{\beta}_{\alpha}p(x)dx})\leq\frac{\int^{\beta}_{\alpha}\phi(k(x))p(x)dx}{\int^{\beta}_{\alpha}p(x)dx},
\end{eqnarray*}
whenever the right-hand side exists and is finite. Equality occurs only when $k(x)\equiv \,constant$.
\end{lemma}

\section{Proof of Theorem \ref{the2.1}}\setcounter{equation}{0}
We divide the proof into 2 parts.

\textbf{Part 1.} Without the winding number condition $deg(q_i-q_j)\neq0 \,(i\neq j)$, we prove the existence of non-collision solution that forms an Euler's central configuration at any instant.

Note that in the new loop space $\Lambda_2$,
if we let the periodic and non-collision solution of fixed energy system
\begin{eqnarray}\label{4.1}
\cases{m_i\ddot{q}_i+\sum_{{j\neq i\atop{{1\leq j\leq 3}}}}\frac{m_im_j}{|q_j-q_i|^{2}}\frac{q_i-q_j}{|q_j-q_i|}=0,  \quad q_{i}\in\mathbb{R}^{2}, \quad \,\,i=1,2,3,\cr
\frac{1}{2}\sum_{i=1}^{3} m_i|\dot{q}_i(t)|^{2}-\sum_{1\leq i< j\leq3}\frac{m_im_j}{|q_i-q_j|}=h,\,\,\,\,h\in(-\frac{b(\lambda_{0})}{2},\,0)}
\end{eqnarray}
be $\tilde{q}$, then $\tilde{q}$ is also the periodic and non-collision solution of system (\ref{1.1}), i.e. system
\begin{eqnarray*}
m_i\ddot{q}_i+\sum_{1\leq i< j\leq3}\frac{m_im_j}{|q_j-q_i|^{2}}\frac{q_i-q_j}{|q_j-q_i|}=0, \quad q_{i}\in\mathbb{R}^{2}, \quad \,\,i=1,2,3.
\end{eqnarray*}
Employing Remark \ref{rem2.2}, we know the periodic solution and non-collision $\tilde{q}$, is just an Euler's central configuration at any instant. Thus next, we need to prove the existence of $\tilde{q}$.

In fact, in the loop space $\Lambda_2$, since $\|q\|=\[\int_{0}^{T}\(\sum_{i=1}^{3}m_{i}|\dot{q}_{i}|^{2}\)dt\]^{1/2}$, then it is nature to define the functional
\begin{eqnarray*}
\varphi_{1}(q)=\frac{1}{2}\|q\|^{2}\int_{0}^{T}\(h-V(q)\)dt=\frac{1}{2}\|q\|^{2}\int_{0}^{T}\(h+\sum_{1\leq i< j\leq3}\frac{m_{i}m_{j}}{|q_{i}-q_{j}|}\)dt,
\end{eqnarray*}
because the existence of critical points of $\varphi_{1}$ implies the existence of periodic and non-collision
solutions of system (\ref{4.1}). If we take $\tilde{f}(q)=\varphi_{1}(q)$ in Remark \ref{rem3.2}, then in what follows, we only need to find the critical point of $\varphi_{1}$ in $\Lambda_2$. Obviously, our new loop space $\Lambda_2$ is not complete, which implies that if $\tilde{f}(q)$ satisfies the $(P.S.)$ condition in $\Lambda_2$, the critical point (i.e., the limit of the $(P.S.)$ sequence) of $\tilde{f}$ may not belong to $\Lambda_2$, so we can not take $\tilde{f}(q)=\varphi_{1}(q)$ directly. In order to overcome this problem,
according to the perturbation of Newtonian potential $V(q)$,
we substitute
\begin{eqnarray*}
V_{\tilde{\varepsilon}}(q)=-\sum_{1\leq i< j\leq3}\frac{m_{i}m_{j}}{|q_{i}-q_{j}|}+\frac{\tilde{\varepsilon}}{h}\cdot\sum_{1\leq i< j\leq3}\frac{m_{i}m_{j}}{|q_{i}-q_{j}|^{2}}, \quad \tilde{\varepsilon}>0
\end{eqnarray*}
for
\begin{eqnarray*}
V(q)=\sum_{1\leq i< j\leq3}\frac{m_{i}m_{j}}{|q_{i}-q_{j}|},
\end{eqnarray*}
and then consider the functional $\tilde{f}(q)=\varphi_{\tilde{\varepsilon}}(q)$ by (\ref{3.3}). By Lemmas \ref{lem3.3} and \ref{lem3.8}, similar to proof of Theorem A \cite{Ambrosetti1992} (Pages 197-198), we prove that for all $-s/2<h<0$, system (\ref{4.1}) has a periodic and non-collision solution $\tilde{q}\in\Lambda_{2}$ with the value of functional $\varphi_{\tilde{\varepsilon}}(\tilde{q})=\hat{c}=\inf_{\gamma\in\hat{\Gamma}}\max_{t\in[0, 1]}\varphi_{\tilde{\varepsilon}}\big(\gamma(t)\big)$ where
\begin{eqnarray*}
\hat{\Gamma}:=\{\gamma\in C\big([0, 1], \Lambda_2\big):\gamma(0)=q^{(\theta)}, \gamma(\frac{1}{2})=q^{(e_1)},\gamma(1)=q^{(e)}\}.
\end{eqnarray*}
Then from the periodic and non-collision solution $\tilde{q}$ of system (\ref{1.1}), we know that $\tilde{q}$ is also an Euler's central configuration at any instant, and we also do not need the winding number condition $deg(q_i-q_j)\neq0 \,(i\neq j)$ in the new loop space $\Lambda_{2}$.

\textbf{Part 2.} We prove that the non-collision solution $\tilde{q}$ obtained in \textbf{Part 1}, is different from the minimizer of $f(q)$ on $\overline{\Lambda_1}$.

Let the minimizer of $f(q)$ on $\overline{\Lambda_1}$ be $\bar{q}$. Next, we prove $\tilde{q}\neq\bar{q}$ by the contradiction argument, and we assume that
$\tilde{q}=\bar{q}$.

Since $\bar{q}=(\bar{q}_1,\,\bar{q}_2,\,\bar{q}_3)$ is the the periodic and non-collision solution obtained in \cite{Zhang2004}, then by Lemma \ref{lem3.11}, we know that $f(\bar{q})$ is the minimum of the action functional $f(q)$ on $\overline{\Lambda_1}$. Combining (\ref{3.10}), (\ref{3.11}) and Lemma \ref{lem3.10}, we know the minimum of $f_1(q)$ equals to
$$
f_1(\bar{q})=\frac32\cdot(2\pi)^{\frac23}\[b(\lambda_0)/a(\lambda_0)\]^{\frac23}T^{\frac13}.
$$

From \textbf{Part 1} in the proof of Theorem \ref{the2.1}, we know that
$\tilde{q}\in \Lambda_2$. From $q=(q_1,\,q_2,\,q_3)$ and $q(t+\frac{T}{2})=-q(t)$ in $\Lambda_2$, we have $\int^{T}_{0}[q_2(t)-q_1(t)]dt=0$. Therefore
by (\ref{3.11}) and Lemma \ref{lem3.12}, we see that for $q\in\Lambda_2$,
\begin{eqnarray}\label{4.2}
f_{1}(q)=\int_{0}^{T}\Big[\frac{1}{2}|\dot{q}_{1}-\dot{q}_{2}|^{2}
+\frac{b(\lambda_0)/a(\lambda_0)}{|q_1-q_2|}\Big]dt
\geq\frac{3}{2}\cdot(2\pi)^{\frac{2}{3}}\[b(\lambda_0)/a(\lambda_0)\]^{\frac{2}{3}}T^{\frac{1}{3}}.
\end{eqnarray}

On the other hand, by Lemma \ref{lem3.12} in \cite[Lemma 9.1.2]{Zhang2018} (Lines 5-6, page 158) and Lemma \ref{lem3.13}, we can claim that if (\ref{4.2}) takes equality, then $|q_2(t)-q_{1}(t)|\equiv \,constant$.
To prove this claim, we take $\alpha=0,\, \beta=T,\,p(t)\equiv1, \,k(t)=|q_1-q_2|^{2},\,\phi(y)=y^{-1/2}, \,y=|q_1-q_2|^{2}\,(q_1\neq q_2)$ in Lemma \ref{lem3.13}, then
\begin{eqnarray}\label{4.3}
\frac{\int_{0}^{T}\frac{1}{|q_1-q_2|}dt}{T}\geq\(\frac{\int_{0}^{T}|q_1-q_2|^{2}dt}{T}\)^{-\frac{1}{2}},
\end{eqnarray}
and thus
\begin{eqnarray}\label{4.4}
\int_{0}^{T}\frac{b(\lambda_0)/a(\lambda_0)}{|q_1-q_2|}dt\geq\frac{b(\lambda_0)}{a(\lambda_0)}\(\int_{0}^{T}|q_1-q_2|^{2}dt\)^{-\frac{1}{2}}T^{\frac{3}{2}}.
\end{eqnarray}
From $\int^{T}_{0}[q_2(t)-q_1(t)]dt=0$, we have
\begin{eqnarray}\label{4.5}
\int_{0}^{T}\frac{1}{2}|\dot{q}_{1}-\dot{q}_{2}|^{2}dt\geq\frac{2\pi^{2}}{T^{2}}\int_{0}^{T}|q_{1}-q_{2}|^{2}dt.
\end{eqnarray}
By (\ref{3.11}), (\ref{4.4}) and (\ref{4.5}), we have
\begin{eqnarray}\label{4.6}
f_{1}(q)&=&\int_{0}^{T}\Big[\frac{1}{2}|\dot{q}_{1}-\dot{q}_{2}|^{2}
+\frac{b(\lambda_0)/a(\lambda_0)}{|q_1-q_2|}\Big]dt \nonumber\\
&\geq&\frac{2\pi^{2}}{T^{2}}\int_{0}^{T}|q_{1}-q_{2}|^{2}dt+\frac{b(\lambda_0)}{a(\lambda_0)}\(\int_{0}^{T}|q_1-q_2|^{2}dt\)^{-\frac{1}{2}}T^{\frac{3}{2}}.
\end{eqnarray}
Let $S=(\int_{0}^{T}|q_{1}-q_{2}|^{2}dt)^{1/2}$, from (\ref{4.6}), one computes that
\begin{eqnarray}\label{4.7}
f_{1}(q)\geq\frac{b(\lambda_0)}{a(\lambda_0)}T^{\frac{3}{2}}\cdot S^{-1}+\frac{2\pi^{2}}{T^{2}}\cdot S^{2}\triangleq g(S).
\end{eqnarray}
By a direct computation, we have
\begin{eqnarray}\label{4.8}
\min \,{g(S)}=g(S_0)=\frac{3}{2}\cdot(2\pi)^{\frac{2}{3}}\[\frac{b(\lambda_0)}{a(\lambda_0)}\]^{\frac{2}{3}}T^{\frac{1}{3}}, \ \ \ with \ \ \
S_0=(2\pi)^{-\frac{2}{3}}\Big[\frac{b(\lambda_0)}{a(\lambda_0)}\Big]^{\frac{1}{3}}T^{\frac{7}{6}}.
\end{eqnarray}
From (\ref{4.3})-(\ref{4.8}), (\ref{4.2}) takes equality if and only if (\ref{4.4}) takes equality and $S=S_0$. On the other hand, with the aid of Lemma \ref{lem3.13}, (\ref{4.4}) takes equality if and only if $|q_2(t)-q_{1}(t)|\equiv \,constant$. Hence, if (\ref{4.2}) takes equality, then $|q_2(t)-q_{1}(t)|\equiv \,constant$. Thus there is no loss of generality in assuming that $|q_2(t)-q_{1}(t)|=l$. We also note that in this situation,
\begin{eqnarray*}
S_0=(\int_{0}^{T}|q_{1}-q_{2}|^{2}dt)^{1/2}=(2\pi)^{-\frac{2}{3}}\Big[\frac{b(\lambda_0)}{a(\lambda_0)}\Big]^{\frac{1}{3}}T^{\frac{7}{6}}.
\end{eqnarray*}
So
\begin{eqnarray}\label{4.9}
|q_2(t)-q_{1}(t)|\equiv l=(2\pi)^{-\frac{2}{3}}\Big[\frac{b(\lambda_0)}{a(\lambda_0)}\Big]^{\frac{1}{3}}T^{\frac{2}{3}}.
\end{eqnarray}

By $\tilde{q}=\bar{q}$, (\ref{4.4}), (\ref{4.5}) and the above analysis, then we know that
\begin{eqnarray}\label{4.10}
\cases{\int_{0}^{T}\frac{b(\lambda_0)/a(\lambda_0)}{|q_1-q_2|}dt=\frac{b(\lambda_0)}{a(\lambda_0)}\(\int_{0}^{T}|q_1-q_2|^{2}dt\)^{-\frac{1}{2}}T^{\frac{3}{2}},\cr
\int_{0}^{T}\frac{1}{2}|\dot{q}_{1}-\dot{q}_{2}|^{2}dt=\frac{2\pi^{2}}{T^{2}}\int_{0}^{T}|q_{1}-q_{2}|^{2}dt.}
\end{eqnarray}

Noticing that our solution $\tilde{q}$ of system (\ref{1.1}), is obtained by using the extended mountain pass theorem for fixed energy system (\ref{4.1}), so the total energy $h\in(-b(\lambda_{0})/2,\,0)$, then
combining (\ref{2.2}), (\ref{3.9})-(\ref{3.10}), (\ref{4.9})-(\ref{4.10}) and
\begin{eqnarray*}
h=\frac{1}{2}\sum_{i=1}^{3} m_i|\dot{q}_i(t)|^{2}-\sum_{1\leq i< j\leq3}\frac{m_im_j}{|q_i-q_j|},
\end{eqnarray*}
we have
\begin{eqnarray*}
-\frac{b(\lambda_0)}{2}T<\int_{0}^{T}hdt&=&a(\lambda_0)\int_{0}^{T}\Big[\frac{1}{2}|\dot{q}_{1}
-\dot{q}_{2}|^{2}-\frac{b(\lambda_0)/a(\lambda_0)}{|q_1-q_2|}\Big]dt\nonumber\\
&=&a(\lambda_0)\frac{2\pi^{2}}{T^{2}}\int_{0}^{T}|q_{1}-q_{2}|^{2}dt}-b(\lambda_0)\(\int_{0}^{T}
|q_1-q_2|^{2}dt\)^{-\frac{1}{2}}T^{\frac{3}{2}\nonumber\\
&=&a(\lambda_0)(2\pi^{2})(2\pi)^{-\frac{4}{3}}T^{\frac{1}{3}}\Big[\frac{b(\lambda_0)}{a(\lambda_0)}
\Big]^{\frac{2}{3}}
-b(\lambda_0)(2\pi)^{\frac{2}{3}}T^{\frac{1}{3}}\Big[\frac{b(\lambda_0)}{a(\lambda_0)}\Big]^{-\frac{1}{3}}
\nonumber\\
&=&-b(\lambda_0)\Big[\frac{a(\lambda_0)}{2b(\lambda_0)}\Big]^{\frac{1}{3}}\pi^{\frac{2}{3}}
T^{\frac{1}{3}}\nonumber\\
&<&-b(\lambda_0)T,
\end{eqnarray*}
and it is a contradiction.

Thus, (\ref{4.2}) can not take equality, which implies $\tilde{q}=\bar{q}$ is impossible. Therefore the solution $\tilde{q}$ obtained in \textbf{Part 1},
is different from the variational
minimizer $\bar{q}$ of the Lagrangian action on the loop space $\overline{\Lambda_1}$.  $\Box$

\begin{remark}\label{rem4.1}
In \cite{Zhang2004}, by using a direct variational method, the authors proved that system (\ref{1.1}) has a periodic solution $\bar{q}$ in $\overline{\Lambda_1}$ that three bodies form Euler's central configuration at any instant, but now employing the study of fixed energy system (\ref{4.1}) and the extended mountain pass theorem (i.e. Lemma \ref{lem3.1}), we can prove that for all $-b(\lambda_0)/2<h<0$, where $b(\lambda_0)$ is only depended on the masses of the three bodies, system (\ref{1.1}) has an periodic and non-collision solution $\tilde{q}$ in $\Lambda_{2}$ that three bodies form Euler's central configuration at any instant. Moreover, combining the suitable periodic $T$, we prove that the new solution $\tilde{q}$ in $\Lambda_{2}$ that three bodies form Euler's central configuration at any instant,
is not the variational
minimizer of the Lagrangian action on the loop space $\overline{\Lambda_1}$.
\end{remark}

\vskip3mm\noindent
\textbf{Acknowledgements}
\vskip2mm\noindent
The first author is partially supported by NSF of China (11671278, 12071316) and the
research funding project of Guizhou Minzu University (GZMU[2019]QN04). The third author is partially supported by NSF of China (11501577).

\end{document}